\documentclass{article}
\usepackage{amsmath,amsfonts,amsthm}
\newtheorem{thm}{Theorem}[section]

\newtheorem{cor}[thm]{Corollary}

\begin{document}

\title{\bf Lifting Lie Algebras over the Residue Field of a Discrete
Valuation Ring}

\maketitle
\centerline{William Browder$^{a,}$\footnote{Supported by the National 
Science Foundation.} \& Jonathan 
Pakianathan$^{b,}$\footnote{Supported by the Sloan Doctoral Dissertation
Fellowship.}}

\begin{align*}
\begin{split}
a. \ Dept.\ of\ Math.,\ Princeton\ University,\ Princeton,\ NJ\ 08544. \\
b. \ Dept.\ of\ Math.,\ University\ of\ Wisconsin,\ Madison,\ WI\ 53706. 
\end{split}
\end{align*}

\begin{abstract}
Let $\mathbf{R}$ be a discrete valuation ring with maximal ideal 
$\pi\mathbf{R}$ and residue field $\mathbf{k}$.
We study obstructions to lifting a Lie algebra $\mathfrak{L}$ over 
$\mathbf{R}/\pi^k\mathbf{R}$ to
one over $\mathbf{R}/\pi^{k+1}\mathbf{R}$. If $\bar{\mathfrak{L}}$ is the
Lie algebra over $\mathbf{k}$ obtained from reducing $\mathfrak{L}$ then
we show there exists a well-defined class in $H^3(\bar{\mathfrak{L}},ad)$
which vanishes if and only if $\mathfrak{L}$ lifts. Furthermore, if $\mathfrak{L}$ lifts, the lifts are shown to be in one to one correspondence with
$H^2(\bar{\mathfrak{L}},ad)$.

\end{abstract}

\section{Introduction.}

Let $\mathbf{R}$ be a Discrete Valuation Ring (DVR) with maximal ideal
$\pi\mathbf{R}$. Let $\mathbf{k}$ be its residue field.
For convenience let $\mathbf{R}_k=\mathbf{R}/\pi^k\mathbf{R}$ for $k \geq 1$. 
An example to keep in mind is $\mathbf{R}=\hat{\mathbb{Z}}_p$, the p-adic integers.
In this case, $\mathbf{R}_k=\mathbb{Z}_{p^k}$, the ring of integers modulo
$p^k$ and the residue field
is $\mathbb{F}_p$, the field with $p$ elements.

A Lie algebra $\mathfrak{L}$ over $\mathbf{R}_k$, $k \geq 1$ is a free
 $\mathbf{R}_k$-module of finite
rank, equipped with a bilinear map 
$$
[\cdot,\cdot] : \mathfrak{L}
 \otimes \mathfrak{L} \rightarrow
\mathfrak{L}
$$ 
which is alternating, i.e., $[v,v]=0$ for all $v \in \mathfrak{L}$, and which 
satisfies the Jacobi identity:
$$
[[u,v],w] + [[v,w],u] + [[w,u],v] =0
$$
for all $u,v,w \in \mathfrak{L}$. 

Given such a Lie algebra we denote
$\bar{\mathfrak{L}}$ the $\mathbf{k}$-Lie algebra obtained by setting 
$$
\bar{\mathfrak{L}}=\mathfrak{L} \underset{\mathbf{R}}{\otimes} \mathbf{k}.
$$
Let $n=dim(\bar{\mathfrak{L}})$.
 Recall the existance of a differential Koszul complex:
$$
0 \rightarrow \wedge^0(\bar{\mathfrak{L}},ad) \overset{d}{\rightarrow} \dots
\overset{d}{\rightarrow}
 \wedge^n(\bar{\mathfrak{L}}, ad)
\rightarrow 0
$$
where $\wedge^i(\bar{\mathfrak{L}},ad)$ is the vector space of alternating
 $\bar{\mathfrak{L}}$-valued i-forms on $\bar{\mathfrak{L}}$, and d is given by
the formula
\begin{align*}
\begin{split}
(d\omega)(x_0,\dots, x_s)= &\sum_{i=0}^s(-1)^i[x_i,\omega(x_0,\dots,\Hat{x_i},
\dots,x_s)] + \\ 
&  \sum_{i<j}(-1)^{i+j}\omega([x_i,x_j],x_0,\dots,\Hat{x_i},\dots,\Hat{x_j},
\dots,x_s)
\end{split}
\end{align*}
for a s-form $\omega$. The cohomology of this complex is the Lie algebra
cohomology $H^*(\bar{\mathfrak{L}},ad)$. 

In this paper we will study the question of when there exists
a $\mathbf{R}_{k+1}$-Lie algebra $\Hat{\mathfrak{L}}$ lifting $\mathfrak{L}$,
 i.e.,
$$
\Hat{\mathfrak{L}} \overset{mod}{\rightarrow} \mathfrak{L}
$$ 
is a surjective map which preserves the brackets. We will say if this 
happens that $\mathfrak{L}$ lifts to a Lie algebra
over $\mathbf{R}_{k+1}$.

This question is answered in 
theorem~\ref{thm: obstruction} of section~\ref{sec: obstruction}. 
Suppose $\bar{\mathfrak{L}}$ is the reduction of
$\mathfrak{L}$ to a Lie algebra over the residue field $\mathbf{k}$.
Then we will show that there exists a class in $H^3(\bar{\mathfrak{L}},ad)$
such that $\mathfrak{L}$ lifts if and only if this class vanishes.

If $\mathfrak{L}$ lifts, one can ask if the lift is unique and if not,
how many distinct ones there are. This question is answered in 
theorem~\ref{thm: lifts} of section~\ref{sec: countinglifts}.
We show that if $\mathfrak{L}$ lifts then the distinct lifts are in one to one
correspondence with $H^2(\bar{\mathfrak{L}},ad)$. 

Finally in section~\ref{sec:example} we provide an example of a Lie
algebra which does not lift.

\section{The obstruction to lifting.}
\label{sec: obstruction}
Let $\mathfrak{L}$ be a Lie algebra over $\mathbf{R}_k$.

There obviously exists a free
$\mathbf{R}_{k+1}$-module $B$ with
$B \overset{mod}{\rightarrow} \mathfrak{L}$ surjective. Let $\{ \bar{e_1},\dots,
\bar{e_n} \}$ be a basis for this free $\mathbf{R}_{k+1}$-module $B$.
 Then
if we set 
$$
e_i=mod(\bar{e_i}),
$$ 
then $\{ e_1, \dots, e_n \}$ is a $\mathbf{R}_k$
basis for $\mathfrak{L}$. 

So we can define a bilinear map
$ [\cdot,\cdot] : B \otimes B \rightarrow B$ by defining on the basis elements:
\begin{align*}
\begin{split}
[\bar{e_i},\bar{e_j}] &= \text{a lift of } [e_i,e_j] \text{ for } i<j \\
[\bar{e_i},\bar{e_i}] &= 0 \\
[\bar{e_i},\bar{e_j}] &= -[\bar{e_j},\bar{e_i}] \text{ for } i>j.
\end{split}
\end{align*}
It is now routine to show that 
$$
[\cdot,\cdot] : B \otimes B \rightarrow B
$$
is alternating, i.e., $[b,b]=0$ for all $b \in B$, and also that 
$$
B \overset{mod}{\rightarrow} \mathfrak{L}
$$ 
preserves the brackets.
We will call such a pair $(B, [\cdot,\cdot])$ a bracket algebra 
(Lie algebra minus Jacobi identity) lift of $\mathfrak{L}$.

If we have another bracket $[\cdot,\cdot]'$ on $B$
lifting the bracket of $\mathfrak{L}$,
then 
$$
[\cdot,\cdot]'-[\cdot,\cdot]=\pi^k\langle \cdot,\cdot\rangle 
$$ 
where 
$$
\pi^k\langle \cdot,\cdot\rangle : B \otimes B \rightarrow B
$$ 
is an alternating, bilinear
form which takes values in $\pi^kB$. 

As $\pi^{k+1}B = 0$, we see that
$$
\pi^k\langle \pi x,y\rangle=0
$$ 
for all $x,y \in B$. 
Define 
$$
\chi : \pi^kB \rightarrow \bar{\mathfrak{L}}
$$ 
by $\chi(\pi^kx)=\lambda(x)$ where 
$$
\lambda : B \rightarrow \bar{\mathfrak{L}}
$$
is the mod reduction from $\mathbf{R}_{k+1}$ to $\mathbf{k}$. It is easy to
check that $\chi$ is a well-defined isomorphism of abelian groups. 
Let $\psi$ be its inverse.
So we see easily that
$\chi(\pi^k\langle \cdot,\cdot\rangle)$ induces an alternating, 
bilinear form
$$
\langle \cdot,\cdot\rangle : \bar{\mathfrak{L}} \otimes \bar{\mathfrak{L}} \rightarrow
 \bar{\mathfrak{L}}.
$$

Conversely, given such an alternating bilinear
form 
$$
\langle \cdot,\cdot\rangle : \bar{\mathfrak{L}} \otimes \bar{\mathfrak{L}}
\rightarrow \bar{\mathfrak{L}},
$$ 
we can set 
$$
[\cdot,\cdot]'=[\cdot,\cdot]
+ \psi(\langle \lambda(\cdot),\lambda(\cdot)\rangle)
$$ 
to get another lift of $\mathfrak{L}$.


Now suppose we have a bracket
algebra lift $(B, [\cdot,\cdot])$ of $\mathfrak{L}$.
Then we can define
$$
\bar{J}_{[\cdot,\cdot]}(x,y,z)=[[x,y],z]+[[y,z],x]+[[z,x],y]
$$
for all $x,y,z \in B$.

It is easy to check that $\bar{J}_{[\cdot,\cdot]}$
defines a $B$-valued alternating 3-form on $B$. Moreover we note that
$B$ lifts $\mathfrak{L}$ which is a Lie algebra hence we see
$\bar{J}_{[\cdot,\cdot]}$ has image in $\pi^kB$. Hence using $\chi$
as mentioned before we see that
$\chi(\bar{J}_{[\cdot,\cdot]})$ induces an $\bar{\mathfrak{L}}$-valued
alternating 3-form on $\bar{\mathfrak{L}}$ which we will call
$J_{[\cdot,\cdot]}$.

It is easy to see that
$J_{[\cdot,\cdot]}=0$ if and only if $\bar{J}_{[\cdot,\cdot]}=0$ if and only
if $(B,[\cdot,\cdot])$ is a Lie algebra.

Note $J_{[\cdot,\cdot]}$ is an element of $\wedge^3(\bar{\mathfrak{L}},ad)$
of the Koszul complex mentioned above. Let us show it is closed under
the differential of the Koszul complex. For simplicity we drop the subscript
and denote $J_{[\cdot,\cdot]}$ by $J$ and similarly $\bar{J}_{[\cdot,\cdot]}$
by $\bar{J}$. Then
\begin{align*}
\begin{split}
(dJ)(x,y,z,w)= &- J([x,y],z,w) + J([x,z],y,w) - J([x,w],y,z) \\
	      &- J([y,z],x,w) + J([y,w],x,z) - J([z,w],x,y) \\
	     &+ [x,J(y,z,w)] - [y,J(x,z,w)] \\
	     &+ [z,J(x,y,w)] - [w,J(x,y,z)] 
\end{split}
\end{align*}
for all $x,y,z \in \bar{\mathfrak{L}}$.

Now we apply $\psi$ to both sides and note 
$$
\psi(J(x,y,z))=\bar{J}(\bar{x},\bar{y},\bar{z})
$$ 
where $\bar{x},\bar{y},\bar{z}$ are lifts of $x,y,z$ to $B$, so we get:
\begin{align*}
\begin{split}
\psi(dJ(x,y,z,w)) =& -\bar{J}([\bar{x},\bar{y}],\bar{z},\bar{w})
+ \bar{J}([\bar{x},\bar{z}],\bar{y},\bar{w})
 - \bar{J}([\bar{x},\bar{w}],\bar{y},\bar{z}) \\
& -\bar{J}([\bar{y},\bar{z}],\bar{x},\bar{w})
+ \bar{J}([\bar{y},\bar{w}],\bar{x},\bar{z})
- \bar{J}([\bar{z},\bar{w}],\bar{x},\bar{y}) \\
& + [\bar{x},\bar{J}(\bar{y},\bar{z},\bar{w})] 
- [\bar{y},\bar{J}(\bar{x},\bar{z},\bar{w})] \\
& + [\bar{z},\bar{J}(\bar{x},\bar{y},\bar{w})]
- [\bar{w},\bar{J}(\bar{x},\bar{y},\bar{z})] \\
=&  0
\end{split}
\end{align*}
where to get the last equality we plugged in the definition of
$\bar{J}$ and cancelled terms in pairs
using that $[\cdot,\cdot]$ is an alternating form. 

Thus we conclude that $dJ=0$, i.e., 
$J_{[\cdot,\cdot]}$ is closed under the differential 
of the Koszul complex. 

Given another bracket $[\cdot,\cdot]'$ on $B$
which lifts $\mathfrak{L}$ we want to relate $J_{[\cdot,\cdot]}$ and
$J_{[\cdot,\cdot]'}$. Recall 
$$
[\cdot,\cdot]'=[\cdot,\cdot]
+ \psi(\langle \lambda(\cdot),\lambda(\cdot)\rangle)
$$ 
for some alternating, bilinear two form
$$
\langle \cdot,\cdot\rangle : \bar{\mathfrak{L}} \otimes \bar{\mathfrak{L}} \rightarrow
\bar{\mathfrak{L}}
$$ 
which we can view as an element in 
$\wedge^2(\bar{\mathfrak{L}},ad)$ of the Koszul complex. Then we have:
\begin{align*}
\begin{split}
\psi(J_{[\cdot,\cdot]'}(x,y,z)) =&
\bar{J}_{[\cdot,\cdot]'}(\bar{x},\bar{y},\bar{z}) \\
=& [[\bar{x},\bar{y}]',\bar{z}]' + [[\bar{y},\bar{z}]',\bar{x}]'
+ [[\bar{z},\bar{x}]',\bar{y}]' \\
=& [[\bar{x},\bar{y}],\bar{z}]' + [\psi(\langle x,y\rangle),\bar{z}]' \\
& [[\bar{y},\bar{z}],\bar{x}]' + [\psi(\langle y,z\rangle),\bar{x}]' \\
& [[\bar{z},\bar{x}],\bar{y}]' + [\psi(\langle z,x\rangle),\bar{y}]' \\
=& [[\bar{x},\bar{y}],\bar{z}] + \psi(\langle [x,y],z\rangle)
+ [\psi(\langle x,y\rangle),\bar{z}] \\
& + \psi(\langle \lambda(\psi(\langle x,y\rangle)),z\rangle) \\
& +[[\bar{y},\bar{z}],\bar{x}] + \psi(\langle [y,z],x\rangle)
+ [\psi(\langle y,z\rangle),\bar{x}] \\
& + \psi(\langle \lambda(\psi(\langle y,z\rangle)),x\rangle) \\
& +[[\bar{z},\bar{x}],\bar{y}] + \psi(\langle [z,x],y\rangle)
 + [\psi(\langle z,x\rangle),\bar{y}] \\
& + \psi(\langle \lambda(\psi(\langle z,x\rangle)),y\rangle) 
\end{split}
\end{align*}
where $\lambda : B \rightarrow \bar{\mathfrak{L}}$ is the mod map.
 Since $k \geq 1$,
we see easily that $\lambda \circ \psi = 0$ so the equation above
 simplifies
to:
\begin{align*}
\begin{split}
\psi(J_{[\cdot,\cdot]'}(x,y,z)) = &
\bar{J}_{[\cdot,\cdot]}(\bar{x},\bar{y},\bar{z}) \\
& + \psi(\langle [x,y],z\rangle + \langle [y,z],x\rangle + \langle [z,x],y\rangle) \\
& + [\psi(\langle x,y\rangle),\bar{z}] + [\psi(\langle y,z\rangle),\bar{x}] \\
& + [\psi(\langle z,x\rangle),\bar{y}] \\
= & \psi(J_{[\cdot,\cdot]}(x,y,z)) \\
& + \psi(\langle [x,y],z\rangle + \langle [y,z],x\rangle + \langle [z,x],y\rangle) \\
& + \psi([\langle x,y\rangle,z] + [\langle y,z\rangle,x] + [\langle z,x\rangle,y]) \\
= & \psi(J_{[\cdot,\cdot]}(x,y,z) - (d\langle \cdot,\cdot\rangle)(x,y,z))
\end{split}
\end{align*}
where in the second to last equality we used the easily verified identity
$$ 
[\psi(\langle x,y\rangle),\bar{z}]=\psi([\langle x,y\rangle,z])
$$ 
for all $x,y,z \in \bar{\mathfrak{L}}$ and in the last equality we thought 
of $\langle \cdot,\cdot\rangle$
as an element in $\wedge^2(\bar{\mathfrak{L}},ad)$ of the Koszul complex.

We apply $\chi$ to both sides of the equation above to get that
$$
J_{[\cdot,\cdot]'}=J_{[\cdot,\cdot]} -d(\langle \cdot,\cdot\rangle)
$$ 
so $J_{[\cdot,\cdot]'}$ and $J_{[\cdot,\cdot]}$ differ by a boundary
in the Koszul complex and hence represent the same cohomology class $[J]$
in $H^3(\bar{\mathfrak{L}},ad)$.
The following theorem follows immediately:

\begin{thm}
\label{thm: obstruction}
Given a Lie algebra $\mathfrak{L}$ over $\mathbf{R}_k$, let $\bar{\mathfrak{L}}$
be its reduction to a $\mathbf{k}$-Lie algebra as before. Then there
exists a well-defined cohomology class 
$$
[J] \in H^3(\bar{\mathfrak{L}},ad)
$$ 
which vanishes if and only if $\mathfrak{L}$ has a $\mathbf{R}_{k+1}$-Lie 
algebra lift $\Hat{L}$ in the sense mentioned before. 
\end{thm}

A p-adic Lie-algebra $\mathfrak{L}$ will be a finitely generated,
free $\Hat{\mathbb{Z}}_p$-module equipped with a bracket as usual which satisfies
the Jacobi identity,
where $\Hat{\mathbb{Z}}_p$ is the p-adic integers. Such a $\mathfrak{L}$ reduces
in an obvious way to $\bar{\mathfrak{L}}$, a Lie algebra
over $\mathbb{F}_p=\mathbb{Z}_p$. Conversely given $\bar{\mathfrak{L}}$ a Lie
algebra over $\mathbb{F}_p$ we say it lifts to the p-adics if such
a $\mathfrak{L}$ exists. Now we can state the following corollary as an
immediate consequence of the theorem above.

\begin{cor}
A Lie algebra $\bar{\mathfrak{L}}$ over $\mathbb{F}_p$ with
$H^3(\bar{\mathfrak{L}},ad)=0$ lifts to the p-adics.
\end{cor}

The obstruction given in theorem~\ref{thm: obstruction} is nontrivial
in general as shown by the example in section~\ref{sec:example}.
However if the reduction map 
$\mathbf{R}_{k+1} \overset{mod}{\longrightarrow} \mathbf{R}_k$ splits,
that is, if there exists a ring map 
$$
\mathbf{R}_k \overset{\kappa}{\longrightarrow} \mathbf{R}_{k+1}
$$ 
with $mod \circ \kappa = Identity$, then Lie algebras over $\mathbf{R}_k$ always lift to ones over $\mathbf{R}_{k+1}$ using $\kappa$ so the obstructions obtained must all vanish.
This happens for example when $\mathbf{R}=\mathbb{F}_p[[x]]$ the power
series ring over $\mathbb{F}_p$. Note this ring has residue field
$\mathbb{F}_p$ the same residue field as for the p-adic integers
where this obstruction theory is nontrivial (see section~\ref{sec:example}).

\section{Classifying lifts of $\mathfrak{L}$.}
\label{sec: countinglifts}

Suppose for this section that we are in the situation as in the
last section but that $\mathfrak{L}$ has a Lie algebra lift $\Hat{\mathfrak{L}}_0$
to a $\mathbf{R}_{k+1}$-Lie algebra. 

Suppose now that $\Hat{\mathfrak{L}}_1$,
$\Hat{\mathfrak{L}}_2$ are two other lifts. We will view these two lifts as 
having the same underlying
$\mathbf{R}_{k+1}$-module $B$ but with different bracket structures
$[\cdot,\cdot]_1$ and $[\cdot,\cdot]_2$ respectively. 

Then as we saw before,
there are two alternating 2-forms $\langle \cdot,\cdot\rangle_i  : \bar{\mathfrak{L}} \otimes
\bar{\mathfrak{L}} \rightarrow \bar{\mathfrak{L}}$ such that 
$$
[\cdot,\cdot]_i=[\cdot,\cdot]_0 + 
\psi(\langle \lambda(\cdot),\lambda(\cdot)\rangle_i)
$$
for $i = 1,2$.

Now as 
$$
J_{[\cdot,\cdot]_i} = J_{[\cdot,\cdot]_0} - d(\langle \cdot,\cdot\rangle_i)
$$
and 
$$
J_{[\cdot,\cdot]_0}=J_{[\cdot,\cdot]_1}=J_{[\cdot,\cdot]_2}=0,
$$ 
we see that 
$\langle \cdot,\cdot\rangle_i$ are closed elements of $\wedge^2(\bar{\mathfrak{L}},ad)$ for $i=1,2$.

Suppose that they determine the same cohomology element in 
$H^2(\bar{\mathfrak{L}},ad)$, i.e., 
$$
\langle \cdot,\cdot\rangle_2-\langle \cdot,\cdot\rangle_1 = -d\phi
$$ 
where
$\phi \in \wedge^1(\bar{\mathfrak{L}},ad)$ is a linear map $\bar{\mathfrak{L}}
\rightarrow \bar{\mathfrak{L}}$. Then 
$$
\Psi = Id + \psi \circ \phi \circ \lambda
 : B \rightarrow B
$$ 
is certainly a homomorphism of $\mathbf{R}_{k+1}$-
modules. (Here Id stands for the identity map). On the other hand it is
invertible using 
$$
(1+x)^{-1} = \sum_{i=0}^{\infty}(-x)^i
$$ 
as 
$\psi \circ \phi \circ \lambda$ takes values in $\pi^kB$ and hence
has square zero. So $\Psi$ is an automorphism of $B$. 

Let us show
$[\Psi(x),\Psi(y)]_2 = \Psi([x,y]_1)$ so that $\Psi$ gives an isomorphism
between the Lie algebras $\Hat{\mathfrak{L}}_1$ and $\Hat{\mathfrak{L}}_2$.
\begin{align*}
\begin{split}
[\Psi(x),\Psi(y)]_2 &= [x + \psi(\phi(\lambda(x))), y +
 \psi(\phi(\lambda(y)))]_2 \\
&= [x,y]_2 + [x, \psi(\phi(\lambda(y)))]_2 +
 [\psi(\phi(\lambda(x))),y]_2 
\end{split}
\end{align*}
where we do not write the 4th term as it is zero because
the image of $\psi$ is in $\pi^kB$.
\begin{align*}
\begin{split} 
[\Psi(x),\Psi(y)]_2 =& [x,y]_0 + \psi(\langle \lambda(x),\lambda(y)\rangle_2) \\
& + [x,\psi(\phi(\lambda(y)))]_0 + \psi(\langle \lambda(x),0\rangle_2) \\
& + [\psi(\phi(\lambda(x))),y]_0 + \psi(\langle 0,\lambda(y)\rangle_2) \\
=& [x,y]_0 + \psi(\langle \lambda(x),\lambda(y)\rangle_2) \\
& + [x,\psi(\phi(\lambda(y)))]_0  + [\psi(\phi(\lambda(x))),y]_0 \\
=& [x,y]_0 + \psi(\langle \lambda(x),\lambda(y)\rangle_2) \\
& + \psi([\lambda(x),\phi(\lambda(y))])
+ \psi([\phi(\lambda(x)),\lambda(y)]).
\end{split}
\end{align*}
On the other hand we have:
\begin{align*}
\begin{split}
\Psi([x,y]_1) =& [x,y]_1 + \psi(\phi(\lambda([x,y]_1))) \\
	      =& [x,y]_1 + \psi(\phi([\lambda(x),\lambda(y)])) \\
	      =& [x,y]_0 + \psi(\langle \lambda(x),\lambda(y)\rangle_1) \\
               & + \psi(\phi([\lambda(x),\lambda(y)])). 	
\end{split}
\end{align*}
So we have:
\begin{align*}
\begin{split}
[\Psi(x),\Psi(y)]_2 - \Psi([x,y]_1) =&
\psi(\langle \lambda(x),\lambda(y)\rangle_2 - \langle \lambda(x),\lambda(y)\rangle_1) + \\
&\psi([\lambda(x),\phi(\lambda(y))] + [\phi(\lambda(x)),\lambda(y)]
-\phi([\lambda(x),\lambda(y)])) \\
=& \psi(\langle \lambda(x),\lambda(y)\rangle_2 - \langle \lambda(x),\lambda(y)\rangle_1) + \\
&\psi((d\phi)(\lambda(x),\lambda(y))),
\end{split}
\end{align*}
where we have used the formula:
$$ (d\phi)(u,v) = -\phi([u,v]) + [u,\phi(v)] - [v,\phi(u)]. $$
Note that up to now to get the formula above we have not used the fact that
$$
\langle \cdot,\cdot\rangle_2 - \langle \cdot,\cdot\rangle_1 = -d\phi
$$ 
and have only used the fact that
$$
\Psi = Id + \psi \circ \phi \circ \lambda.
$$
Now since
$$
\langle \cdot,\cdot\rangle_2 = \langle \cdot,\cdot\rangle_1 - d\phi(\cdot,\cdot),
$$
we see from the above that $[\Psi(x),\Psi(y)]_2 = \Psi([x,y]_1)$.

So we see indeed that $\Hat{\mathfrak{L}}_1$ and $\Hat{\mathfrak{L}}_2$ are isomorphic
Lie algebras under an isomorphism which induces the identity map on
$\mathfrak{L}$.

Conversely suppose $\Hat{\mathfrak{L}}_1$ and $\Hat{\mathfrak{L}}_2$
are lifts which are isomorphic via an isomorphism $\Psi$ which induces
the identity
map on $\mathfrak{L}$ then clearly $\Psi = Id + \mu$ where $\mu$ is
a homomorphism of $\mathbf{R}_{k+1}$-modules which has image
in $\pi^kB$. Thus easy to see 
$$
\Psi = Id + \psi \circ \phi \circ
\lambda
$$ 
for some $\phi : \bar{\mathfrak{L}} \rightarrow \bar{\mathfrak{L}}$
a linear map. Reversing the steps above we see that 
$$
\langle \cdot,\cdot\rangle_2 = \langle \cdot,\cdot\rangle_1 - d\phi
$$ 
so that the elements
in $\wedge^2(\bar{\mathfrak{L}},ad)$ corresponding to the two Lie algebras
$\Hat{\mathfrak{L}}_1$ and $\Hat{\mathfrak{L}}_2$ are cohomologous. So from
all the preceding facts, the following theorem follows readily:

\begin{thm}
\label{thm: lifts}
Given $\mathfrak{L}$ a Lie algebra over $\mathbf{R}_k$, let
$\bar{\mathfrak{L}}$ be the Lie algebra obtained by reducing to the residue
field. Suppose there is a Lie algebra lift of $\mathfrak{L}$ to
a $\mathbf{R}_{k+1}$ Lie algebra. Then the set of all lifts (up to
isomorphism of Lie algebras inducing identity on $\mathfrak{L}$) are
in one to one correspondence with $H^2(\bar{\mathfrak{L}},ad)$.
\end{thm}

\section{A Lie algebra which does not lift.}
\label{sec:example}
Consider the Lie algebra $\mathfrak{M}=\mathfrak{sl_3(3)}$ of trace zero , 3x3 matrices
over $\mathbb{F}_3$, the field with three
 elements. This is an eight dimensional
Lie algebra which is easily verified to be perfect, i.e.,
$[\mathfrak{M},\mathfrak{M}]=\mathfrak{M}$.

Unlike the characteristic zero case, this
Lie algebra has a one dimensional center consisting of the multiples
of the identity which have trace zero. Form the quotient Lie algebra of
$\mathfrak{M}$ by its center and let us call it $\mathfrak{P}$. 

This Lie algebra $\mathfrak{P}$ is seven dimensional and
is still perfect of course. It is easy to check it is in fact simple.
A tedious calculation reveals its Killing form to be identically zero.
(This is much different than in the characteristic zero case where
the Killing form is nondegenerate if and only if the Lie algebra
is semisimple). On the
other hand there does exist a nondegenerate, invariant, symmetric
3-form and this generates the vector space of invariant, symmetric
3-forms on $\mathfrak{P}$. One can also show that $\mathfrak{P}$ is unimodular
and hence calculate its cohomology relatively easily as 
\begin{align*}
\begin{split}
H^1(\mathfrak{P},\mathbb{F}_3) &= 0 \\
H^2(\mathfrak{P},\mathbb{F}_3) &= \mathbb{F}_3^6 \\
H^3(\mathfrak{P},\mathbb{F}_3) &= 0, 
\end{split}
\end{align*} 
and the rest filled in by Poincar$\acute{e}$ duality.
Also one can show 
\begin{align*}
\begin{split}
H^0(\mathfrak{P},ad) &= 0 \\
H^1(\mathfrak{P},ad) &= \mathbb{F}_3^7. 
\end{split}
\end{align*}
Thus $\mathfrak{P}$ is not complete although one can show directly that
it is restricted using the $p$-power map from $\mathfrak{M}$. Furthermore
this is the unique restricted structure.

As is well-known, such cohomology would not be expected for a simple
Lie-algebra in the characteristic zero case. So we have much reason
to believe $\mathfrak{P}$ does not lift to the p-adics. Indeed, computer 
calculations of the cohomology obstruction in section~\ref{sec: obstruction}
have shown that it does not
even lift to a $\mathbb{Z}_9$-Lie algebra. Unfortunately, at the current
time we have left it at that and not attempted to try to find a ``human''
proof. Unfortunately although the verification that $[J]$ is not
a boundary in this case is merely linear algebra, it is of considerable
size!

In general the Lie algebras $\mathfrak{sl}_p(p)$ for $p$ a prime number
have a nontrivial center and the quotient Lie algebra obtained by
moding out this center should be quite interesting.

For an application of the results of this paper to questions about uniform,
p-central, p-groups see \cite{P}.

\section*{Acknowledgements.}

The authors would like to thank the Mathematical Institute of Aarhus
University for its hospitality, and Thomas Weigel for his suggestions.

\end{document}